\documentclass[12pt]{article}

\usepackage{amsmath}
\usepackage{amsfonts}
\usepackage{amssymb}

\newtheorem{theorem}{Theorem}
\newtheorem{prop}{Proposition}

\title{Extension of the interior connection of a nonholonomic
manifold with a Finsler metric}

\author{Sergey V. Galaev}

\begin{document}

\maketitle\vskip-50ex
{\renewcommand{\abstractname}{Abstract}\begin{abstract} The
notions of the interior  and  truncated connections of a
nonholonomic manifold are introduced. A class of extended
truncated connections is distinguished. For the case of a contact
space with a Finsler metric, it is shown that there exists a
unique extended truncated connection that satisfies additional
properties. The curvature tensor of the obtained connection in the
case of a sub-Riemannian space coincides with the Wagner curvature
tensor that was constructed by Wagner for the case of an arbitrary
nonholonomic manifold of codimension one endowed with an interior
affine connection.
\end{abstract}

\section*{Introduction}

A nonholonomic manifold is a smooth distribution on a smooth
manifold. This distribution is in general not integrable. The are
many works devoted to the theory of nonholonomic systems and
nonholonomic geometry, among them there are a lot of works of
prominent geometers of the first half of the twentieth century.
The foundations of nonholonomic geometry were laid in the
classical works on nonholonomic mechanics of physicists and
mathematicians like Hertz, H\"older, Chaplygin, Appell, Korteweg
and others. In the works of Gibbs and Caratheodory on the
foundations of thermodynamics there appears for the first time the
contact nonholonomic structure, which is the simplest nonholonomic
structure. The nonholonomic theory is strongly related with the
theory of distributions that in particular includes the theory of
Pfaffian systems. An important role here plays the works of
\'E.~Cartan, who developed the formalism of exterior differential
forms and distributions. In the twenties of the last century,
after the works of  Levi-Civita and H.~Weyl, where the Riemannian
and affine connections were defined and the relation of geometry
and mechanics was discovered, there appeared the understanding of
what the nonholonomic mechanics should use for finding new
geometric structures. The origin of such interplay was placed by
Vranceanu and Synge. In USSR the nonholonomic theory was actively
developed by V.F\,Kagan. In 1937 he proposed the following topic
for the N.I.\,Lobachevskii prize: "to establish the foundations of
the doctrine for the theory of nonholonomic spaces ... the
applications to mechanics, physics and the theory of integration
of Pfaffian equations are desirable". The most serious
achievements on the nonholonomic geometry and its applications in
mechanics belong to Vranceanu, Synge, Schouten and V.V.\,Wagner.
The Romanian mathematician Vranceanu was the first who precisely
formulated the notion of the nonholonomic structure on a
Riemannian manifold as well as its relations to the dynamics of
nonholonomic systems. Synge studied the problem of the stability
of the inertial motion of nonholonomic mechanical systems and by
this he anticipated the notion of the curvature of a nonholonomic
manifold. Schouten defined the parallel transport of some vectors
along some vector fields, afterwards this was called the truncated
connection. He also defined the curvature tensor of this
connection. The next step in the definition of the curvature
tensor of a nonholonomic manifold belongs to V.V.~Wagner, who
constructed a curvature tensor expanding the Schouten tensor
\cite{W41}. The Wagner curvature tensor is zero if and only if the
Schouten-Vranceanu connection is flat.

At present the increasing interest to the nonholonomic geometry is
mostly related to the active usage of nonholonomic systems
(sub-Riemannian spaces) in the control theory.
The problems of the control theory frequently require the
generalization of the Riemannian metric to the Finslerian one
\cite{H-P}.


\section{Gradient coordinate system on a nonholonomic manifold}

In this paper $X$ denotes a connected $C^\infty$-manifold of
dimension $n=2m+1$, $m\geq 2$. All objects on $X$ are assumed to
be smooth.

{\it A nonholonomic manifold} $D$ is a smooth not involutive
distribution of codimension one on $X$. Assume that there exists a
one-dimensional distribution $D^\bot$ such that
\begin{equation}\label{eq1}TX=D\oplus D^\bot.\end{equation}
Following Wagner, we call $D^\bot$ a closing of the nonholonomic
manifold $D$.

A vector field on $X$ is called  admissible if it is tangent to
the distribution $D$. A 1-form on $X$ is called admissible  if it
is zero on the closing $D^\bot$. At last, an admissible tensor
field on $D$  is a linear combination of tensor products of the
admissible vector fields and 1-forms. Denote by $f^p_q(D)$ the
module of admissible tensor fields of type $(p,q)$ on $D$.

A coordinate chart $K(x^\alpha)$ ($\alpha,\beta,\gamma=1,...,n$;
$a,b,c=1,...,2m$) on the manifold $X$ is called adapted to the
nonholonomic manifold $D$ if $\partial_n=\frac{\partial}{\partial
x^n}\in f^1_0(D^\bot)$. It is not hard to check that any two
adapted coordinate charts are related by a transformation of the
form
\begin{equation}\label{eq2} x^a= x^a (x^{a'}),\qquad x^n= x^n (x^{a'},x^{
n'}).\end{equation} Let $P:TX\to D$ be the projector defined by
the decomposition \eqref{eq1} and let $K(x^\alpha)$ be an adapted
coordinate chart. Then the vector fields
$$P(\partial_a)=\vec{e}_a=\partial_a-\Gamma_a^n\partial_n$$ are
linearly independent and they generate the distribution $D$, i.e.
$D={\rm span} (\vec{e}_a)$. Thus we have on the manifold $X$ the
nonholonomic basis field $(\vec{e}_a,\vec{e}_n)$ and the
corresponding cobasis field $(dx^a,\Theta^n=dx^n+ \Gamma_a^n
dx^a)$. The vector fields $\vec{e}_a$ defined in  the nonholonomic
manifold linear coordinates called by Wagner gradient coordinates.

The gradient coordinates satisfy the property that the
transformation matrix in the formula $\vec{e}_a=A^{ a'}_a\vec{e}_{
a'}$ coincides with the matrix $\frac{\partial x^{a'}}{\partial
x^{a}}$. This property implies in particular the following
proposition.

 \begin{prop} If $t$ is an admissible tensor field of type $(p,q)$, then the object  $\tilde{t}$
 with the components
${\partial }_nt^{a_1\cdots a_p}_{b_1\cdots b_q}$ is an admissible
tensor field of the same type. \end{prop}

{\bf Proof.} Let $t$ be an admissible tensor field of type
$(p,q)$, than in adapted coordinates it holds $$t=t^{a_1\dots a_p
}_{b_1\dots b_q} \vec e_{a_1}\otimes\cdots \otimes \vec
e_{a_p}\otimes dx^{b_1}\otimes \cdots \otimes dx^{b_q}.$$ The
transformation \eqref{eq2} implies
$$t^{ a_1'\dots  a_p' }_{b_1'\dots b_q'}= A^{a'_1}_{a_1}\cdots
A^{a'_p}_{a_q}A^{b_1}_{b'_1}\cdots A^{b_q}_{b'_q}t^{a_1\dots
a_p}_{b_1\dots b_q}.$$ Equalities \eqref{eq2} show that the
functions $\frac{\partial x^{a'}}{\partial x^a}$ do not depend on
the coordinate $x^{n}$. This proves the proposition. $\Box$

\section{Interior and extended connections on a nonholonomic
manifold}

 Developing
the geometry of nonholonomic manifolds, V.V.~Wagner, introduces
the notion of the interior geometry of a nonholonomic manifold $D$
as the collection of the properties of the objects defined on $D$
that depend only on $D$ and its clothing \cite{W41}. The parallel
transport in a nonholonomic manifold is defined by a connection
$\nabla$ that in the terminology of V.V.~Wagner   is called
interior. In some works in addition to the interior connections
one considers also connections for which the parallel transport of
vectors from $D$ along arbitrary curves in $X$ is defined. Such
connections are called connections in the vector bundle $D$
defined be the nonholonomic manifold.  Defining in \cite{W40} the
curvature of the interior connection $\nabla$, V.V.~Wagner
constructs in a special way a connection in the vector bundle $D$.
The curvature tensor of this connection was called afterwards the
Wagner curvature tensor. In this paper we interpret  the Wagner
curvature tensor as the tensor of the nonholonomicity of a smooth
distribution and then, using this interpretation, we define an
analog of the Wagner curvature tensor for nonholonomic manifolds
with Finsler metrics \cite{GCh00}.

 Suppose that on $X$ a contact structure
$\lambda$ is given,  i.e. $\lambda$ is a 1-form such that the rank
of the 2-form $\omega=d\lambda$ equals $2m$, and there exists a
decomposition $$TX=D\oplus D^\bot,$$ where $D=\ker\lambda$ and
$D^\bot=\ker\omega$. We call the nonholonomic manifold $D$ a
contact space.

We say that on the nonholonomic manifold $D$ an interior
connection is given, if the distribution $\tilde{D}={\pi
}^{-1}_*(D)$, where $\pi:D\to X$ is the natural projection, is
decomposed into the direct sum
\begin{equation} \label{eq5}
\tilde{D}=HD\oplus VD,
\end{equation}
where $VD$ is an vertical distribution on the total space $D$.
Thus the interior connection is uniquely  defined by the object
 $G^a_b(x^a,x^n)$ such that  $HD={\rm span}
(\vec{\varepsilon}_{a})$, where
\begin{equation} \label{eq6}
\vec{\varepsilon}_{a}=\partial_a-\Gamma^n_a\partial_n-G^b_a\partial_{n+b}.
\end{equation}
The equality  \eqref{eq6} implies the equalities of the form
\begin{equation} \label{eq7}
\left[\vec{\varepsilon}_a,\vec{\varepsilon}_b\right]={\omega
}_{ba}\partial_n +2(\vec{e}_{[b}G^c_{a]}-G^d_{[a}G^c_{b].d}),
\end{equation}
\begin{equation} \label{eq8}
\left[\vec{\varepsilon}_a,\partial_n\right]=\partial_n\Gamma
^n_a\partial_n+(\partial_nG^b_a)\partial_{n+b}.
\end{equation}
In  \eqref{eq7} a dot denotes the derivative with respect to the
coordinates on a fibre. Denote by $T^p_q(D)$ the set of all
admissible tensors of type $(p,q)$. By definition, an admissible
Finslerian tensor field of type $(p,q)$ is a morphism $t:V \to
T^p_q(D)$ such that $t(z)\in T^p_{\pi (z) q}(D)$, where
$V=D\backslash \{0\}$ is the bundle of non-zero vectors from the
distribution $D$.

The partial case, when the interior connection is given by a
linear connection $(G^a_b)(x^\alpha,x^{n+a})=\Gamma
^a_{bc}(x^\alpha)x^{n+a}$, the objects
$$P^b_a=\partial_nG^b_a,\quad K^c_{ab}=2(\vec{e}_{[b}G^c_{a]}-G^d_{[a}G^c_{b].d})$$ were called
by V.V.~Wagner the first and the second Schouten curvature
tensors, respectively. We use these terms  for the object defined
above.

Let $\vec{u}$ be the vector field defined by the conditions
$\vec{u}\in f_0^1 (D^\bot)$ and $\lambda(\vec{u})=1.$ If we assume
that an adapted coordinate  chart satisfies the additional
condition $\partial_n=\vec{u}$, then the transformation
\eqref{eq2} takes the following more simple form:
\begin{equation}\label{eq3} x^a= x^a (x^{a'}),\qquad x^n= x^{
n'}+const.\end{equation} In what follows we will consider only
adapted coordinate charts related by the transformation
\eqref{eq3}. It is not hard to check the equality
\begin{equation} \label{eq4}
\left[\vec{e}_a,\vec{e}_b \right]=\omega_{ba}\partial_n.
\end{equation}

The distribution $D$ may be considered as the total space of the
vector bundle $\mu =(X,D,\pi )$, where $\pi :D\to X $ is the
natural projection.

Any adapted coordinate chart $k(x^{\alpha })$   on the manifold
$X$ defines the coordinate chart $\tilde{k}(x^{\alpha},x^{n+a})$
on the manifold $\tilde D$, where $x^{n+a}$ are the coordinates of
the vector  $\vec{v}\in D$ with respect to the basis
$(\vec{e}_{a})$. Thus $D$ is a smooth manifold of dimension
$4m+1$.

A coordinate chart on the manifold $D$ defines the nonholonomic
frame field $(\vec{\varepsilon}_a,\partial_n,\partial_{n+a})$,
where $\vec{\varepsilon}_a$ are defined by the equality
\eqref{eq6}. Thus besides the distribution $HD$, which defines the
exterior connection on the nonholonomic manifold $D$, we get the
well defined distribution
\begin{equation}\label{eq9}\widetilde{HD}= HD\oplus{\rm
span}(\partial_n)\end{equation} on the whole manifold $D$. This
distribution defines an infinitesimal connection on $D$ as on the
vector bundle. This means that we deal with a connection that is
traditionally called a truncated connection. We call this
connection  an extension of the interior connection of the
nonholonomic manifold. Any other extension of the interior
connection is defined by a vector field that has the following
coordinate form: $$\vec{u}=\partial_n-G^a_n\partial_{n+a},$$ where
the object $G_n^a$ is an example of an admissible Finslerian
vector field.

\section{Contact space with a Finsler structure}

Suppose now that on the manifold $D$ is defined a function
$L(x^\alpha,x^{n+a})$ that satisfies the following conditions:
\begin{itemize} \item[1)] $L$ is smooth at least on $D\backslash \{0\}$;
\item[2)] $L$  is homogeneous of degree 1 with respect to the coordinates of an admissible vector, i.e.
\begin{equation}\label{eq10}L(x^\alpha,\lambda x^{n+a})=\lambda L(x^\alpha,x^{n+a}),\quad
\lambda>0;\end{equation}
\item[3)] $L(x^\alpha,x^{n+a})$ is positive if not all $x^{n+a}$ are zero simultaneously;
\item[4)] the quadric form $$L^2_{\cdot a\cdot b}\xi^a\xi^b=
\frac{\partial^2L^2}{\partial x^{n+a}\partial x^{n+b}}$$ is
positive definite. \end{itemize} We call the triple $(X,D,F)$,
where $F=L^2$, a contact sub-Finslerian manifold.

The following takes  place.

\begin{theorem} On any contact sub-Finslerian manifold $(X,D,F)$ there exists a unique
trancated metric connection such that
\begin{equation}\label{eq11} G^a_{b\cdot c}=G^a_{c \cdot b}.\end{equation}\end{theorem}

{\bf Proof.} Let $\vec{v}$ be an admissible vector field such that
\begin{align}\label{eq12} \nabla_av^b&=\vec{e}_av^b+G_a^b(x^\alpha,v^c)=0,\\
\label{eq13}
\nabla_nv^b&=\partial_nv^b+G_n^b(x^\alpha,v^c)=0.\end{align} We
need to prove that the objects $G^b_a$, and $G^b_n$ are uniquely
defined by the metrizability condition and by the condition
 \eqref{eq11}. The metrizability condition implies that the function $f$ is constant along
any integrable curve of the vector field $\vec{v}$, hence
$$dF(x^\alpha,x^{n+a})=\vec{e}_aFdx^a+\partial_n F\Theta^n=0,$$
where $x^{n+a}=v^a$. This equality together with \eqref{eq12} and
\eqref{eq13} implies the following:
\begin{align}\label{eq14} \vec{e}_aF-G_a^cF_{\cdot c}&=0,\\
\label{eq15} \vec{e}_nF-G_n^cF_{\cdot c}&=0.\end{align}
Differentiating   the equality \eqref{eq14} with respect to
$x^{n+b}$, we get
$$\vec{e}_aF_{\cdot b}-G^c_{a\cdot b}F_{\cdot c}-G^c_aF_{\cdot c\cdot b}=0.$$
Contracting this with $x^{n+a}$, we obtain
\begin{equation}\label{eq16}
x^{n+a}\vec{e}_aF_{\cdot b}-x^{n+a}G^c_{a\cdot b}F_{\cdot
c}-x^{n+a}G^c_aF_{\cdot c\cdot b}=0.
\end{equation} The homogeneity of the coefficients $G^c_a$ and the condition \eqref{eq11} imply the equality
$$x^{n+a}G^c_{a\cdot b}=G^c_b.$$
Thus \eqref{eq16} takes the form
\begin{equation}\label{eq17}
x^{n+a}\vec{e}_aF_{\cdot b}-x^{n+a}G^c_aF_{\cdot c\cdot
b}-G^c_bF_{\cdot c}=0.\end{equation} On the other hand, from
\eqref{eq14} it follows that $$\vec{e}_aF=G^c_aF_c.$$
Consequently, \eqref{eq17} can be rewritten in the form
\begin{equation}\label{eq18} x^{n+a}G_a^cF_{\cdot c\cdot
b}=x^{n+a}\vec{e}_aF_{\cdot b}-\vec{e}_bF.\end{equation} Define
the following admissible tensor field:
$$g_{ab}=\frac{1}{2}F_{ab}.$$
Then \eqref{eq18} implies
$$x^{n+a}G_a^c=\frac{1}{2}g^{bc}(x^{n+a}\vec{e}_aF_b-\vec{e}_bF).$$
Differentiating this equality by $x^{n+d}$, we get
$$G^c_d+x^{n+a}G^c_{ad}=\frac{1}{2}(g^{bc}(x^{n+a}\vec{e}_aF-\vec{e}_bF))_{\cdot d},$$
or
\begin{equation}\label{eq19} G^c_d=\frac{1}{4}(g^{bc}(x^{n+a}\vec{e}_aF-\vec{e}_bF))_{\cdot
d}.\end{equation} We may find now the components of the field
$G^c_n$ using \eqref{eq15}.

Alternating the second derivative of an admissible vector field
and using \eqref{eq4}, \eqref{eq7}, \eqref{eq17}, we get
\begin{equation}\label{eq20} K^c_{ab}(x^\alpha,v^d)+2\omega_{ab}\partial_nv^c=0.\end{equation}
From this it follows that
\begin{equation}\label{eq21} \partial_nv^d+\omega^{ba}K^d_{ab}=0.\end{equation}
Comparing \eqref{eq21} with \eqref{eq13}, we get
\begin{equation}\label{eq22} G^d_n(x^\alpha,v^c)=\omega^{ba}(x^\alpha)K^d_{ab}(x^\alpha,v^c).\end{equation}

Thus we conclude that if a metric connection satisfying
\eqref{eq11} exists, then its coefficients are uniquely defined by
the equalities \eqref{eq19} and \eqref{eq22}. On the other hand,
defining the coefficients of the connection using \eqref{eq19} and
\eqref{eq22}, we obtain a metric torsion free connection. The
theorem is proved. $\Box$

As it is known, the curvature tensor of the constructed metric
connection on the vector bundle is the tensor of the
nonholonomicity of the corresponding horizontal distribution. In
order to write it down, one finds the Lie brackets of the vector
fields generating the horizontal distribution:

\[[{\vec{\varepsilon }}_a,{\vec{\varepsilon }}_b] = {\omega }_{ba}\vec{U} + R^{c}_{ab}{\partial }_{n+c},\]

\[ [{\vec{\varepsilon }}_a,\vec{U}] = {\partial }_n{\Gamma }^n_a\vec{U} + R^{c}_{na}{\partial }_{n+c},\]
 where

\begin{equation} \label{eqA}
R^{c}_{ab} = K^{c}_{ab} + {\omega }_{ba}{\omega }^{ij}K^{c}_{ij},
\end{equation}

\begin{equation} \label{eqB}
R^{c}_{na} = P^{c}_{na} + \widetilde{\nabla}_aG^c_n.
\end{equation}

The constructed tensor satisfies the usual properties of the
Berwald curvature tensor.

\vskip0.5cm

Saratov State University

Email: sgalaev(at)mail.ru

\end{document}